Dedicated to Lou Kauffman for his 70th birthday

# KNOT THEORY: FROM FOX 3-COLORINGS OF LINKS TO YANG-BAXTER HOMOLOGY AND KHOVANOV HOMOLOGY

## JÓZEF H. PRZYTYCKI

**Abstract:** This paper is an extended account of my "Introductory Plenary talk at Knots in Hellas 2016" conference[1]. We start from the short introduction to Knot Theory from the historical perspective, starting from Heraclas text (the first century AD), mentioning R.Llull (1232-1315), A.Kircher (1602-1680), Leibniz idea of Geometria Situs (1679), and J.B.Listing (student of Gauss) work of 1847. We spend some space on Ralph H. Fox (1913-1973) elementary introduction to diagram colorings (1956). In the second section we describe how Fox work was generalized to distributive colorings (racks and quandles) and eventually in the work of Jones and Turaev to link invariants via Yang-Baxter operators; here the importance of statistical mechanics to topology will be mentioned. Finally we describe recent developments which started with Mikhail Khovanov work on categorification of the Jones polynomial. By analogy to Khovanov homology we build homology of distributive structures (including homology of Fox colorings) and generalize it to homology of Yang-Baxter operators. We speculate, with supporting evidence, on co-cycle invariants of knots coming from Yang-Baxter homology. Here the work of Fenn-Rourke-Sanderson (geometric realization of pre-cubic sets of link diagrams) and Carter-Kamada-Saito (co-cycle invariants of links) will be discussed and expanded. No deep knowledge of Knot Theory, homological algebra, or statistical mechanics is assumed as we work from basic principles. Because of this, some topics will be only briefly described.

## 1. KNOT THEORY STARTED IN PELOPONNESE

As the popular saying goes "All science started in Ancient Greece". Knot Theory is not an exception. We have no proof that ancient Greeks thought of Knot Theory as a part of Mathematics but surgeons for sure thought that knots are important:
a Greek physician named Heraklas, who lived during the first century AD is our main example (see Subsection 1.2). Even before, in pre-Hellenic times, there is mysterious stamp from Lerna, the place famous

---

[1]I decided to keep the original abstract of the talk omitting only the last sentence "But I believe in *Open Talks*, that is I hope to discuss and develop above topics in an after-talk discussion over coffee or tea with willing participants", which applies to a talk but not a paper.

in classical times as the scenes of Herakles' encounter with the hydra [Hea, Wie1, Wie2].

## 1.1. Seal-impressions from Lerna, about 2,200 BC.

Excavations at Lerna by the American School of Classical Studies under the direction of Professor J. L. Caskey (1952-1958) discovered two rich deposits of clay seal-impressions. The second deposit dated from about 2200 BC contains several impressions of knots and links[2] [Hig, Hea, Wie1] (see Fig.1.1).

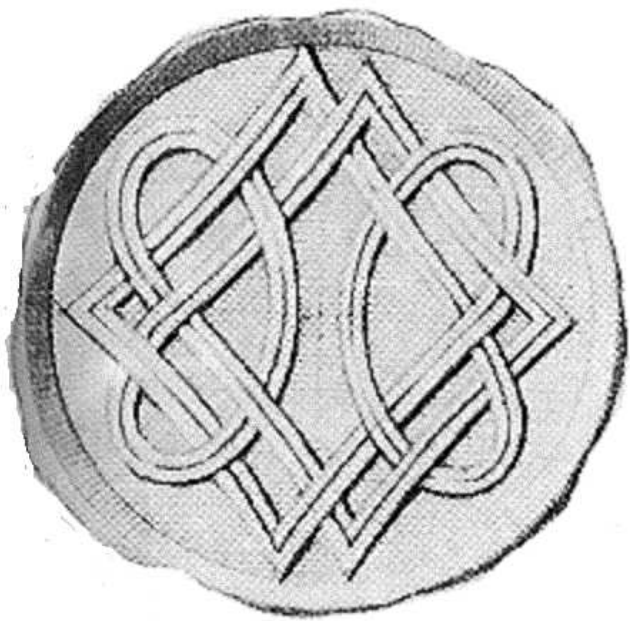

Figure 1.1; A seal-impression from the House of the Tiles in Lerna (c. 2200 BC)

---

[2]The early Bronze Age in Greece is divided, as in Crete and the Cyclades, into three phases. The second phase lasted from 2500 to 2200 BC, and was marked by a considerable increase in prosperity. There were palaces at Lerna, Tiryns, and probably elsewhere, in contact with the Second City of Troy. The end of this phase (in the Peloponnese) was brought about by invasion and mass burnings. The invaders are thought to be the first speakers of the Greek language to arrive in Greece.

I have chosen two more patterns from seals of Lerna; these are not knots or links but "pseudoknots" which I will mention later with respect to extreme Khovanov homology and RNA.

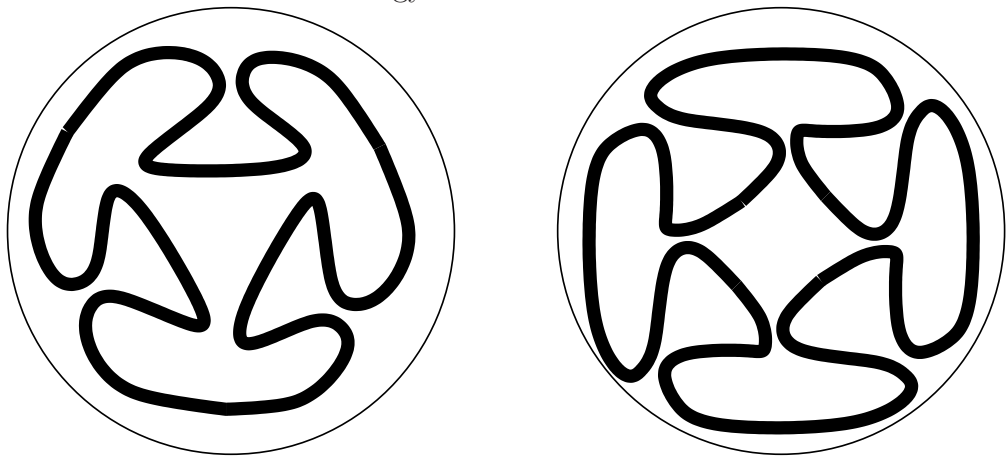

Figure 1.2; pseudoknots from Lerna

## 1.2. **Heraklas slings, first century AD.**

A Greek physician named Heraklas, who lived during the first century AD and who was likely a pupil or associate of Heliodorus[3], wrote an essay on surgeon's slings [He]. Heraklas explains, giving step-by-step instructions, eighteen ways to tie orthopedic slings. Here also Hippocrates "the father of western medicine" should be mentioned[4]. Heraklas work survived because Oribasius of Pergamum (c. 325-400; physician of the emperor Julian the Apostate) included it toward the end of the fourth century in his "Medical Collections".[5] The oldest extant manuscript of "Medical Collections" was made in the tenth century by the Byzantine physician Nicetas. The Codex of Nicetas was brought to Italy in the fifteenth century by an eminent Greek scholar, J. Lascaris, a refugee from Constantinople. Heraklas' part of the Codex of Nicetas has no illustrations[6] and around 1500 an anonymous artist depicted Heraklas' knots in one of the Greek manuscripts of Oribasius "Medical Collections" (in Figure 1.3 we reproduce, after Day and with his comments the first page of drawings [Da-2]). Vidus Vidius (1500-1569), a Florentine who became physician to Francis I (king of France, 1515-1547) and professor of medicine in the Collège de France, translated the Codex of Nicetas into Latin; it contains also drawings of Heraklas' surgeon's slings by the Italian painter, sculptor and architect Francesco Primaticcio (1504-1570); [Da-2, Ra].

---

[3]Heliodorus was a surgeon in the 1st century AD, probably from Egypt, and mentioned in the Satires of Juvenal. This Heliodorus wrote several books on medical technique which have survived in fragments and in the works of Oribasius [Mil]. It is worth to cite Miller: "In the 'Iatrikon Synagogos,' a medical treatise of Oribasius of Pergamum (...) Heliodorus, who lived at the time of Trajan (Roman Emperor 98–117 AD), also mentions in his work knots and loops" [Mil, Sar-1].

[4]Hippocrates of Cos (c. 460-375 BC). A commentary on the Hippocratic treatise on *Joints* was written by Apollonios of Citon (in Cypros), who flourished in Alexandria in the first half of the first century BC. That commentary has obtained a great importance because of an accident in its transmition. A manuscript of it in Florence (Codex Laurentianus) is a Byzantine copy of the ninth century, including surgical illustrations (for example, with reference to reposition methods), which might go back to the time of Apollonios and even Hippocrates. Iconographic tradition of this kind are very rare, because the copying of figures was far more difficult than the writing of the text and was often abandoned [Sar-1]. The story of the illustrations to Apollonios' commentary is described in [Sar-2].

Heraklas' essay is the first surviving text on Knot Theory even if it is not proper Knot Theory but rather its application. The story of the survival of Heraklas' work and efforts to reconstruct his knots in Renaissance is typical of all science disciplines and efforts to recover lost Greek books provided the important engine for development of modern science. This was true in Mathematics as well: the beginning of modern calculus in XVII century can be traced to efforts of reconstructing lost books of Archimedes and other ancient Greek mathematicians. It was only the work of Newton and Leibniz which went much farther than their Greek predecessors.

On a personal note: When I started to be interested in History of Knot Theory, the texts of Heraklas or Oribasius were unknown to Knot Theory community. It was by chance that when in 1992 I had a job interview at Memphis State University I had a meeting with a Dean, an English professor. Learning that I work on Knot Theory he mentioned that he had a friend C.L.Day[7] who wrote a "humanistic" book about knots and their "classical" beginnings. Thus I learned about the book of Day [Da-2] and Heraclas slings.

---

[5]From [Sar-3]:"The purpose of Oribasios *Medical Collection* is so well explained at the beginning of it that it is best to quote his own words": *Autocrator Iulian, I have completed during our stay in Western Gaul the medical summary which your Divinity had commanded me to prepare and which I have drawn exclusively from the writings of Galen. After having praised it, you commanded me to search for and put together all that is most important in the best medical books and all that contributed to attain the medical purpose . I gladly undertook that work being convinced that such a collection would be very useful. (...) As it would be superfluous and even absurd to quote from the authors who have written in the best manner and then again from those who have not written as careful, I shall take my material exclusively from the best authors without omitting anything which I first obtained from Galen....*

[6]Otherwise the Codex of Nicetas is the earliest surviving illustrated surgical codex, containing 30 full page images illustrating the commentary of Appolonios of Kition and 63 smaller images scattered through the pages.

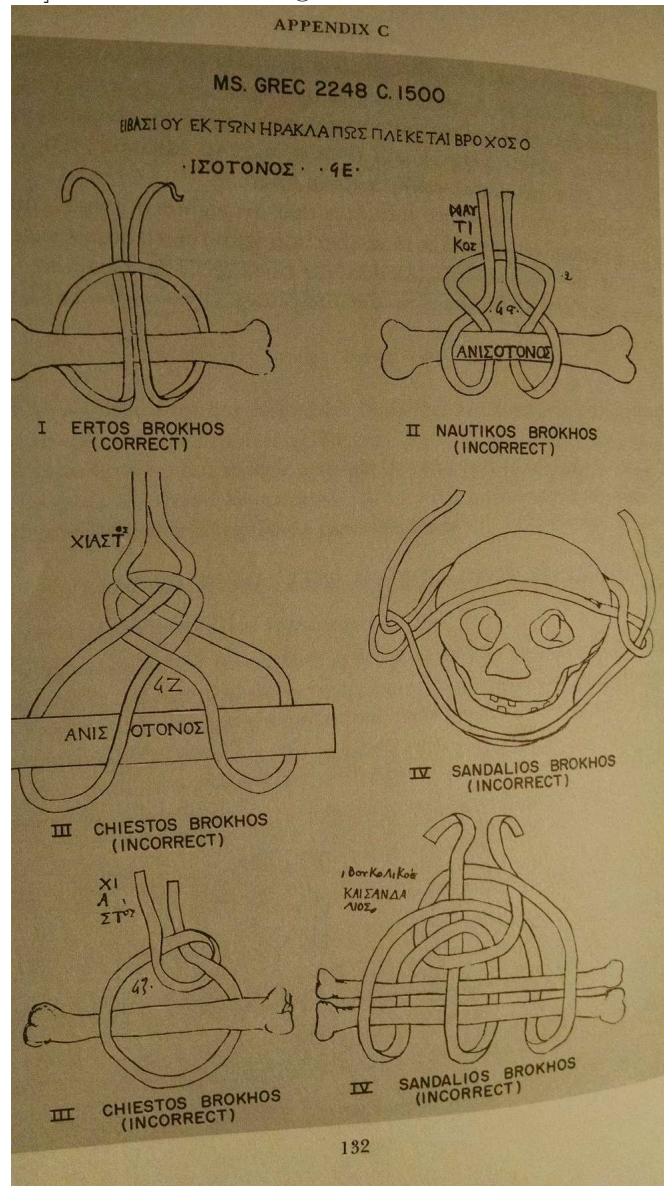

Figure 1.3; slings of Heraklas, c. 100 AD, [Da-2]

### 1.3. **Ramon Llull, Leonardo da Vinci, and Albert Dürer.**

Let us mention in passing the work of Ramon Llull $(1232 - 1315)$ and his combinatorial machines which greatly influenced Leibniz and his idea of Geometria Situs (Figure 1.4).

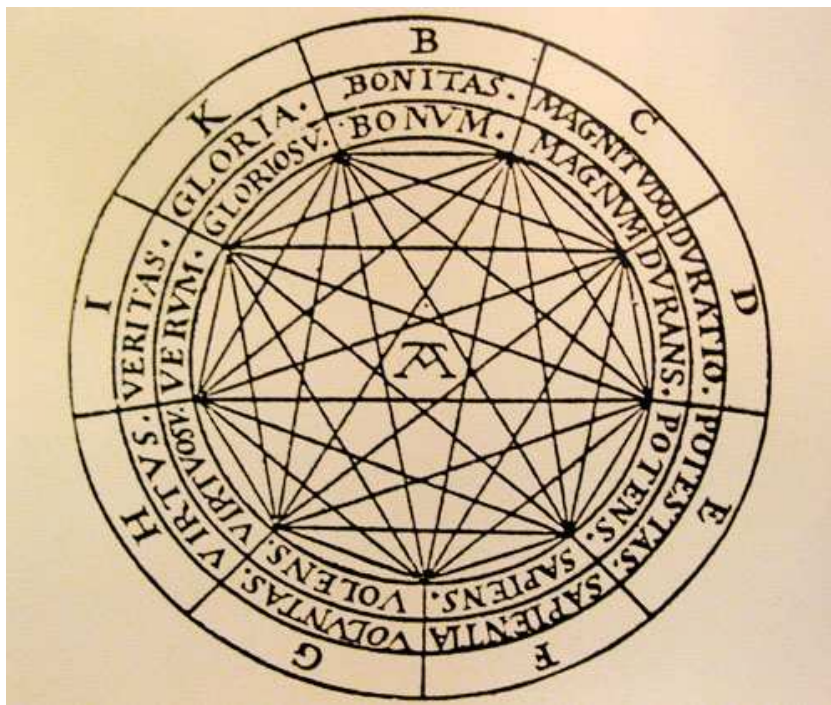

Figure 1.4; Combinatorial machine of Ramon Llull from his Ars Generalis Ultima

The drawing of knots by Leonardo da Vinci and Albert Dürer should be also acknowledged (see Figures 1.5 and 1.6).

Here the anecdote: are this knots really knots or maybe links of more than one component? What is the structure of these links? The graph theorist of note, Frank Harary, took a task of checking it and made precise analysis of Leonardo and Dürer links[8] [Har] (compare also [H-M, H-S]). This paper is based on my Knots in Hellas talks, taking place in ancient Olympia; such a venue is tempting me to write more of history, but I am already straying too far. I would refer to Chapter II of my Book [Prz-8] where I describe the early Knot Theory and the work of Kircher, Leibniz, Vandermonde, Gauss, and Listing; the Chapter is based on my papers [Prz-1, Prz-2].

---

[8]Only after Knots in Hellas conference (July 2016) I learned about the paper by Hoy and Millett [H-M] with very detailed discussion of Leonardo and Dürer knots, see also [H-S].

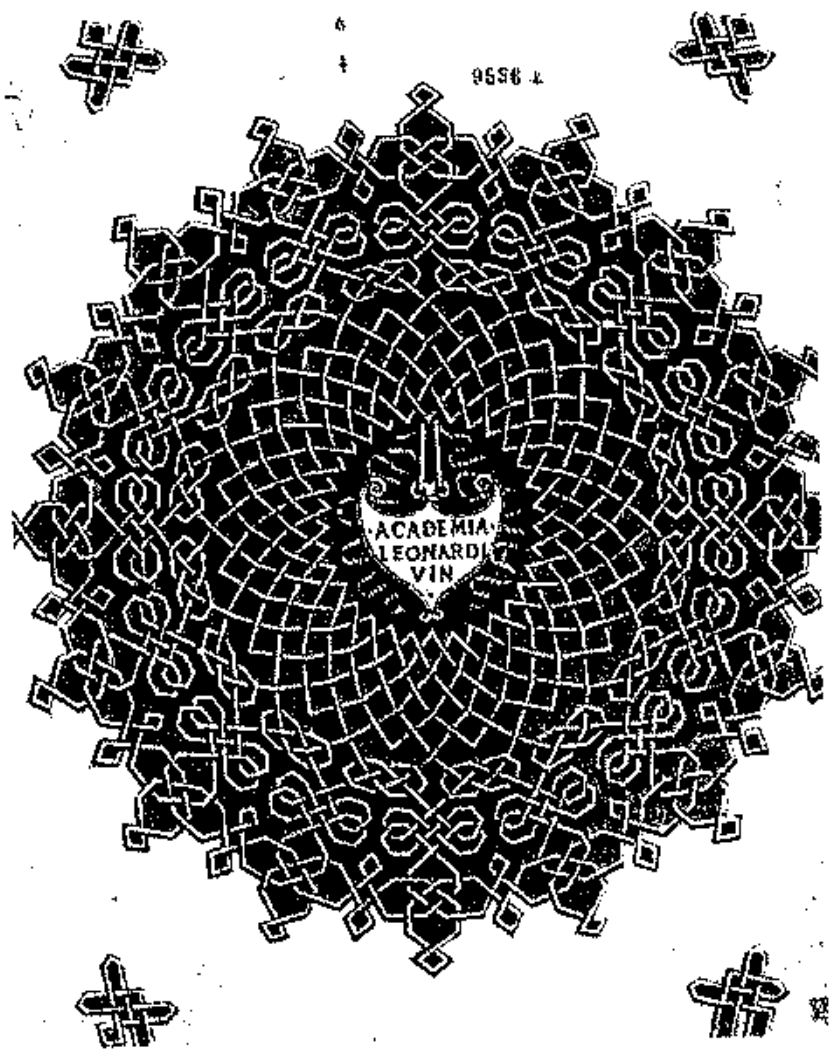

Figure 1.5; Leonardo da Vinci, Leonardus Vinci Academia

Engraving by Leonardo da Vinci[9] (1452-1519) [Mac].

---

[9]Giorgio Vasari writes in [Vas]: "[Leonardo da Vinci] spent much time in making a regular design of a series of knots so that the cord may be traced from one end to the other, the whole filling a round space. There is a fine engraving of this most difficult design, and in the middle are the words: Leonardus Vinci Academia."

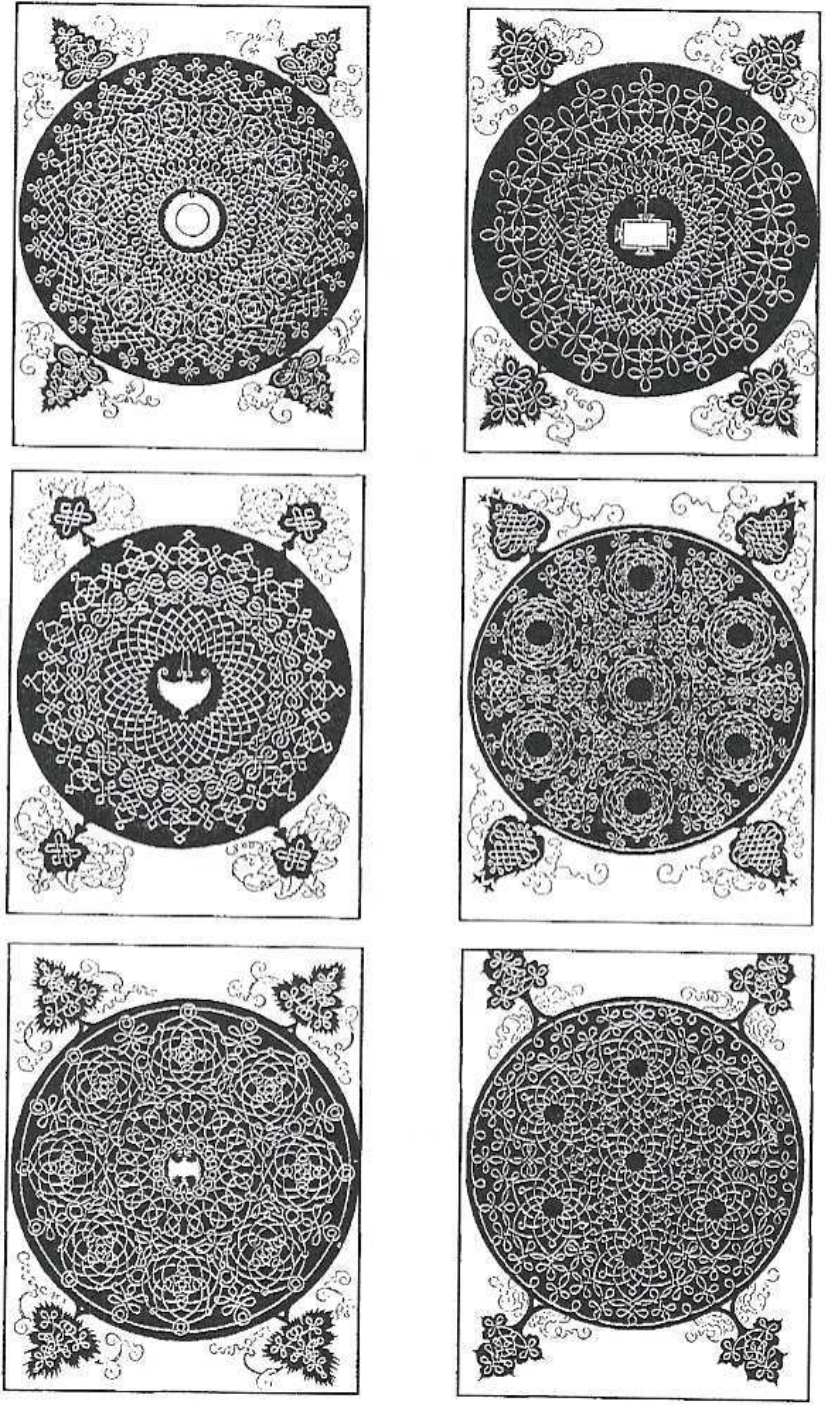

Figure 1.6; Dürer's knots, 1505/6

### 1.4. **Fox and Fox-Trotter colorings.**

My space is finite so let us jump to the second part of XX century, as I promised, to connect Fox colorings with Yang-Baxter invariants (and Khovanov homology).

Yes we all know about Fox colorings; in Figure 1.7 we have iconic nontrivial Fox 3-coloring of the trefoil knot. The rule of 3-coloring is that we color arcs of a diagram using three colors, say red, blue, and yellow in such a way that at each crossing either all colors are used or only one color is used.

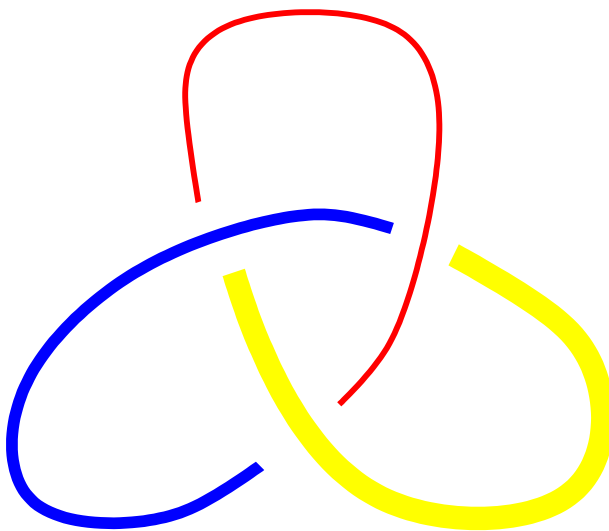

Figure 1.7; nontrivial Fox 3-coloring of the trefoil knot

One also can play the coloring game with the link from Lerna, and yes it has nontrivial Fox 3-coloring as illustrated in Figure 1.8. Indeed Fox 3-colorings motivated many popular, school level articles, notable of them is [Vi-1]. I also wrote about Fox colorings for middle school children [Prz-4, Prz-5].

Still there is some controversy who really invented them[10]. I think it was as follows (I describe likely story based on facts but also my experience with teaching in America). In 1956 Ralph Fox spent a sabbatical at Haveford College as it is explained in the Preface to his book [C-F]: "This book, which is an elaboration of a series of lectures given by Fox at Haveford College while a Philips Visitor there in the spring of 1956, is an attempt to make the subject accessible to everyone. Primarily it is a textbook for a course at the junior-senior level".

---

[10]Reidemeister was considering homomorphisms of the fundamental group of the knot complement into n-dihedral groups. This easily leads to $n$th Fox coloring [Re-1, Re-2].

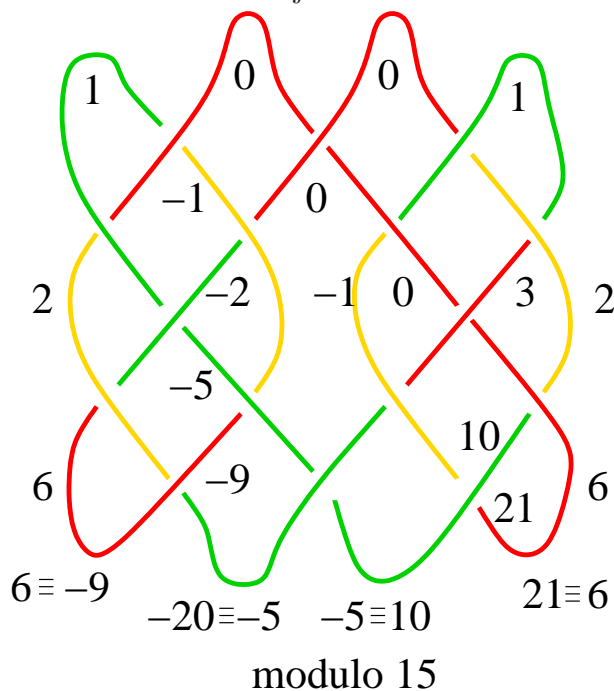

**modulo 15**

Figure 1.8; nontrivial Fox 3-coloring of the Lerna link
Numbers on the picture describe Fox 15-coloring.
In fact the space of colorings of this link is $Col_{\mathbb{Z}}(L) = \mathbb{Z} \oplus \mathbb{Z}_2 \oplus \mathbb{Z}_3 \oplus \mathbb{Z}_5^2 \oplus \mathbb{Z}_7$
where $\mathbb{Z}$ can be represented by trivial (monochromatic) colorings ; compare [Prz-3].

It is curious that 3-coloring and $n$-coloring are mentioned only in Exercises:
For instance Exercises 6 and 7 in Chapter VI are about Fox 3-colorings:
"Exercise 6. Let us say that a knot diagram has property $\ell$ if it is possible to color the projected overpasses in three colors, assigning a color to each edge in such a way that
(a) the three overpasses that meet at a crossing are either all colored the same or are all colored differently;
(b) all three colors are actually used.
Show that a diagram of a knot $K$ has property $\ell$ if and only if $K$ can be mapped homomorphically onto the symmetric group of degree 3.
Exercise 7. Show that property $\ell$ is equivalent to the following: It is possible to assign an integer to each edge in such a way that the sum of the three edges that meet at any crossing is divisible by 3."

Fox 3-coloring can be naturally extended to $n$-coloring, again already hidden in Reidemeister work as homomorphisms of the fundamental group of link complement to the $n$th dihedral group $D_{2n}$, sending meridians to reflections. The diagrammatic definition is hidden in Exercises of Chapter VIII of [C-F], in particular Exercise 8 for $k = n - 1$ and then $-(b - a) = b - c \pmod{n}$ that is around a crossing Fox $n$-coloring looks like $\frac{b}{a} \Big| \frac{c \, \equiv \, 2b-a}{\text{modulo } n}$ .

Fox 11-coloring of the knot $6_2$ of the Rolfsen tables [Rol] is shown in Figure 1.9.

Interestingly when Richard Crowell, a student of Ralph Fox, was talking about Fox colorings to teachers in 1961, he referred to Fox as an inventor of 3-colorings but he said that he learned $n$-coloring from Halle Trotter [Cro]. I asked the question to Trotter and he remembers discussions with Crowell but not inventing $n$-colorings; he kindly answered my inquiry: "I am afraid my historical recollection is now very vague. Dick Crowell and I had many discussions of various things, and if he says so, I perhaps said something to suggest $n$-coloring before he did. He worked out all the details – I did not know even that he was writing the NCTM paper until he sent me a reprint."

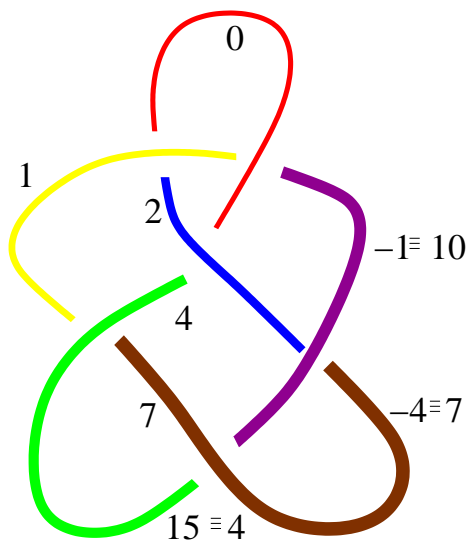

Figure 1.9; nontrivial Fox 11-coloring of the knot $6_2$

## 2. From arc colorings to Yang-Baxter weighted colorings

2.1. **Magma colorings.** Fox $n$-coloring had to wait for its direct generalization (at least in print, discussion of the Conway-Wraith correspondence of 1959 would take another lecture [Wra]) for quarter of century. These, despite the fact that Wirtinger coloring, giving the fundamental group of a link complement or Alexander coloring, giving Alexander module and Alexander polynomial, were known from 1905 and 1928, respectively [Wir, Ale].[11]

With Fox $n$-colorings, an orientation of a diagram is not needed. To save time I will move immediately to oriented diagrams but start naively from a fixed finite set $X$ and coloring arcs of a diagram [12] by elements of $X$. Then we can ask under which conditions the set of such colorings is a link invariant. The simplest approach is now to think that at a crossing colorings change according to some operation. Because we can have positive and negative crossings, we need two operations $*, \bar{*} : X \times X \to X$. The convention for coloring by $(X, *, \bar{*})$ (called 2-magma) is given in Figure 2.1.

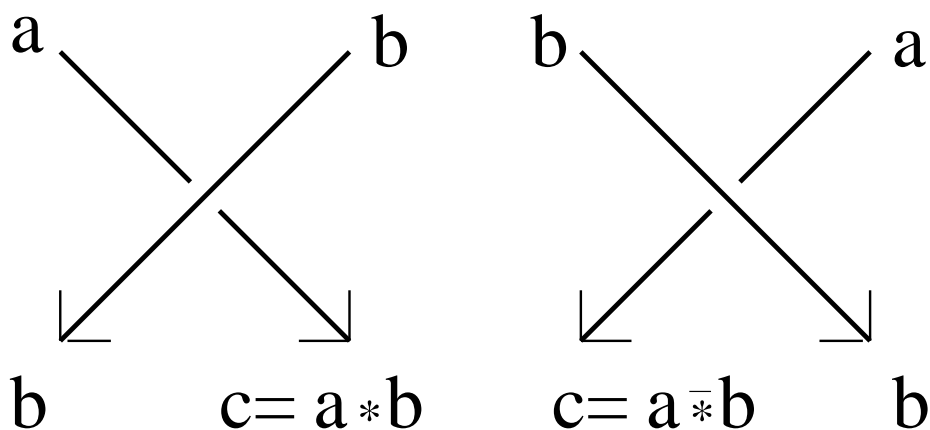

Figure 2.1 Convention for a 2-magma coloring, $f : \mathrm{Arcs}(D) \to X$, of a crossing

[11] The year 1928 is the year of publication of the Alexander's paper, however already in 1919 he discusses in a letter to Oswald Veblen, his former PhD adviser, "a genuine and rather jolly invariant" which we call today the determinant of the knot. It is this construction which Alexander extends later to the Alexander polynomial $\Delta_D(t)$ (determinant is equal to $\Delta_D(t)$ for $t = -1$). In fact the Alexander letter contains more: Alexander constructs the space which we call often today the space of nontrivial Fox $\mathbb{Z}$-colorings or the first homology of the double branched cover of $S^3$ along the knot [A-V].

[12] We consider arcs from undercrossing to undercrossing and semi-arcs from crossing to crossing.

The set of 2-magma colorings is denoted by $Col_X(D)$ and its cardinality by $col_X(D) = |Col_X(D)|$. Of course $col_X(D)$ is not necessary a link invariant and in next subsection we analyze when Reidemeister moves are preserving it.

2.2. **Reidemeister moves and Quandles.** If we want $col_X(D)$ to be a link invariant, we check Reidemeister moves and obtain, after Joyce and Matveev [Joy-2, Mat], the algebraic structure satisfying conditions (1),(2),(3) below, which Joyce in his 1979 PhD thesis named a *quandle* [Joy-1].

**Definition 2.1.**

   (1) $a * a = a$ , for any $a \in X$ *(idempotence condition)*.

   (2) *There is the inverse binary operation*[13] $\bar{*}$, *to* $*$, *that is for any pair* $a, b \in X$ *we have*

$$(a * b)\bar{*}b = a = (a\bar{*}b) * b \ \textit{(invertibility condition)}.$$

      *Equivalently we define* $*_b : X \to X$ *by* $*_b(a) = a * b$, *and invertibility condition means that* $*_b$ *is invertible; we denote* $*_b^{-1}$ *by* $\bar{*}_b$.

   (3) $(a * b) * c = (a * c) * (b * c)$ *(distributivity), for any* $a, b, c \in X$. *Figure 2.5 illustrates how the third Reidemeister move leads to right selfdistributivity, and in fact can be taken as a "proof without words" that* $col_X(D)$ *is preserved by the positive third Reidemeister move if and only if* $*$ *is right self-distributive.*

*If only conditions (2) and (3) hold, then* $(X; *, \bar{*})$ *is called a rack (or wrack); the name coined by J.H.Conway in 1959.*
*If* $* = \bar{*}$ *in the condition (2), that is* $(a * b) * b = a$ *then the quandle is called an involutive quandle or kei* $\equiv\!\!\!\equiv$ *(the last term coined in 1942 by M. Takasaki* [Ta]*).*

Before we show how quandle axioms are motivated by Reidemeister moves it is worth making metamathematical remark:
We have two equivalent approaches to quandle definition. The first approach starts from a magma $(X, *)$ and because the second condition says that $*$ is invertible we can introduce the inverse operation $\bar{*}$. The second approach uses only equations, thus we start from a set $X$ with two binary operations $*$ and $\bar{*}$ and in the second axiom we assume that equations $(a * b)\bar{*}b = a$ and $(a\bar{*}b) * b = b$ hold. In both approaches axioms (1) and (3) are given by equations. The algebraic structure in which conditions are given by identities is called *variety* and G.Birkhoff proved that a class of algebras is a variety if and only if it is closed under homomorphic images, subalgebras, and arbitrary direct products [Bir]. On the other hand the first nonequational approach to a quandle allows that homomorphic image of a quandle is not a quandle (only spindle – the name used for magmas $(X, *)$ satisfying conditions (1) and (3)). See the discussion in Section 1.2 of [Smi] about combinatorial and equational definitions of quasigroup.

---

[13]We can think of "inverse" formally: we introduce the monoid of binary operations on $X$, $Bin(X)$, with composition given by $a(*_1 *_2)b = (a *_1 b) *_2 b$ and identity element $*_0$ given by $a *_0 b = a$, then the inverse means the inverse in the monoid, that is $* \bar{*} = *_0 = \bar{*} *$; see [Prz-6].

After this detour we go back to Reidemeister moves.

**First Reidemeister move and idempotence**

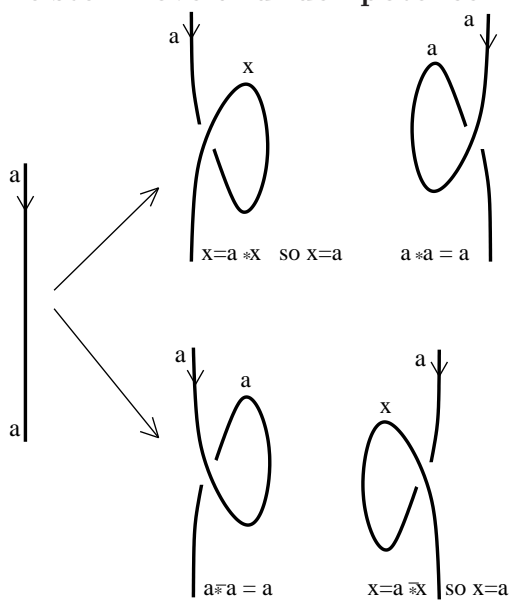

Figure 2.2; First Reidemeister move leads to idempotent conditions $a * a = a$ and $a\bar{*}a = a$
It gives also the stronger condition that $a$ is the only solution of the equation $a * x = x$ and $a\bar{*}x = x$;
However this follows from the idempotent condition and the condition (2) (that $\bar{*}$ is the inverse of $*$)

**First Reidemeister move as framing change**

For many considerations it is important to observe that the first Reidemeister move can be interpreted as a framing change of a framed diagram; Figure 2.3.

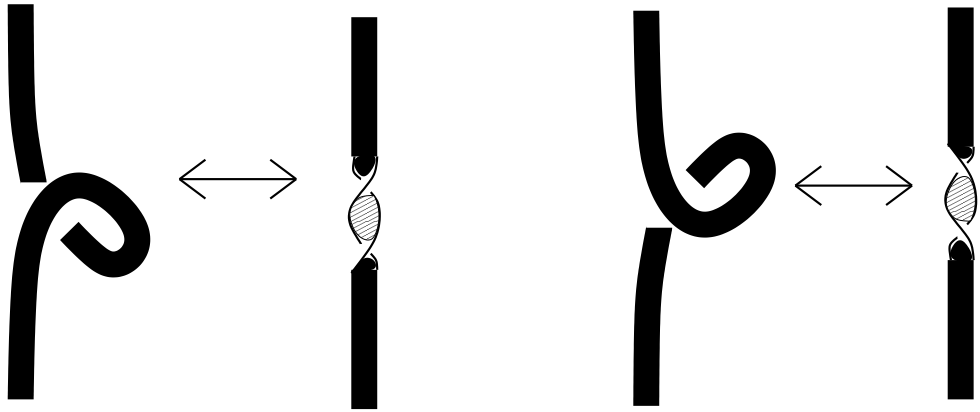

Figure 2.3; first Reidemeister move can be interpreted as a framing change

**Second Reidemeister move and invertibility**

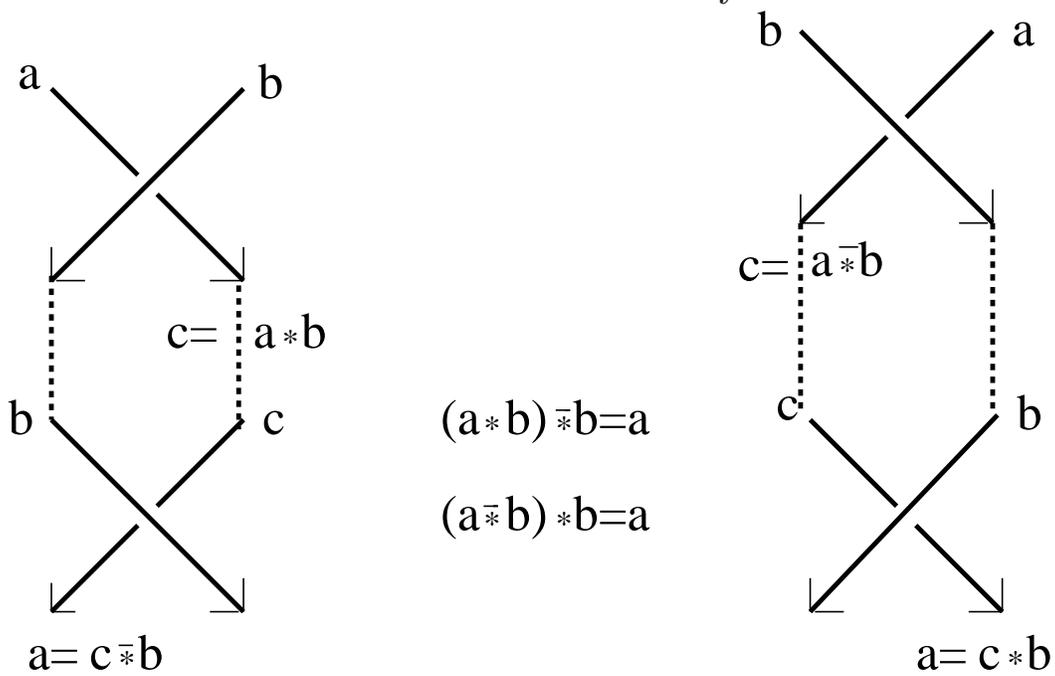

$(a*b)\bar{*}b=a$

$(a\bar{*}b)*b=a$

Figure 2.4; Second Reidemeister move and magma coloring

The cardinality $col_X(D)$ is preserved by the second Reidemeister move if $*$ is invertible.

**Third Reidemeister move and distributivity**

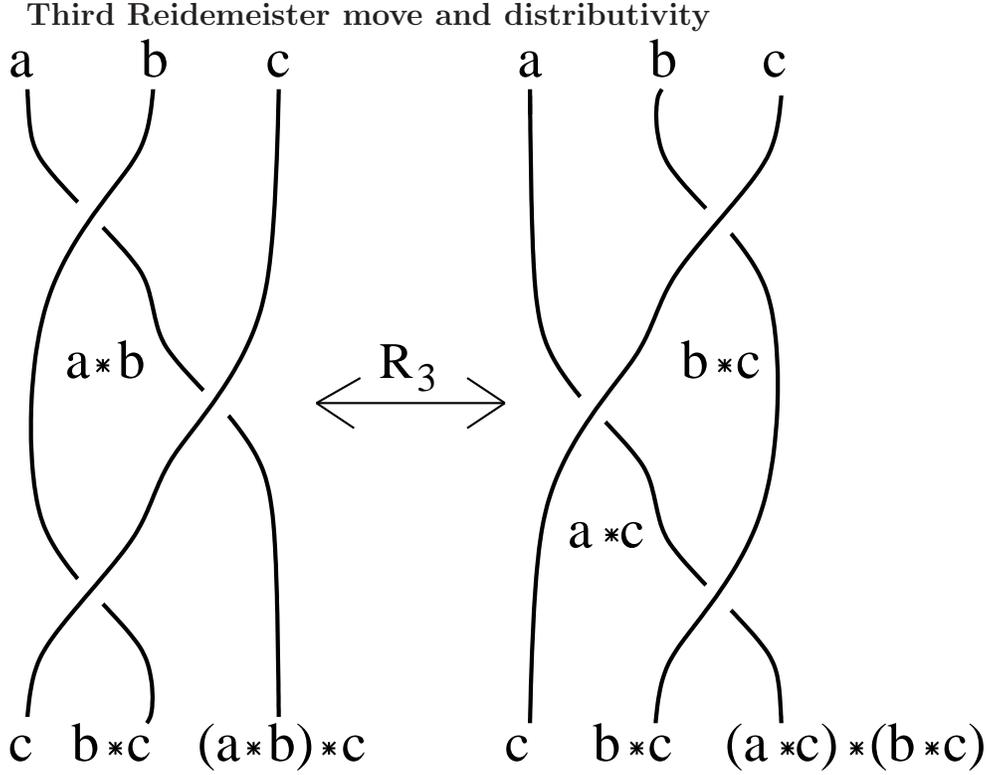

Figure 2.5; Third Reidemeister move leads to right selfdistributivity
$(a * b) * c = (a * c) * (b * c)$

### 2.3. 2-(co)cycle invariants.

Let $X$ be a finite set and $*$ and $\bar{*}$ two binary operations. We define, after Carter-Kamada-Saito [CKS], 2-(co)cycle invariants of links:

(1) A 2-chain, $\Psi(D, \phi)$ associated to the diagram $D$ and coloring of its arcs by $\phi : arcs(D) \to X$ is an element of $\mathbb{Z}X^2$ defined as a sum over all crossings of $D$ of the pair $\pm(a, b)$ according to conditions in figure below, that is $\Psi(D, \phi) = \sum_v sgn(v)(a, b)$, where the sum is taken over all crossings of $D$.

(2) A 2-cochain with coefficients in an abelian group $A$ is a function $\alpha : X^2 \to A$ or equivalently an element of $Hom(\mathbb{Z}X^2, A)$. A 2-cochain associated to the diagram $D$ and coloring of its arcs by $\phi : arcs(D) \to X$ is an element of $Hom(\mathbb{Z}X^2, A)$, defined by $\Psi(D, \phi, \alpha) = \Sigma_v sgn(v)\alpha(a, b)$.

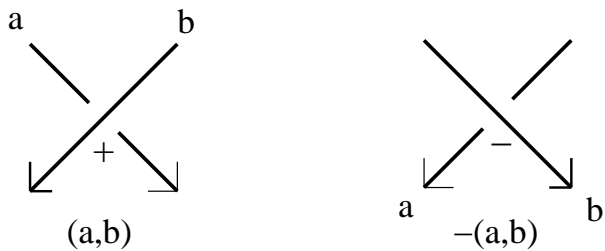

Figure 2.6; Convention for the 2-chain

In order to have Knot Theory applications one would like to have the chain $\Psi(D, \phi)$ (resp. cochain $\alpha$) to be a 2-cycle (resp. cocycle) in some homology (resp. cohomology) theory. Further we would like to have Reidemeister moves preserving homology (resp. cohomology) class. This motivated initially authors of [CJKLS] and led to the discovery that what they need is essentially rack homology introduced around 1990 by Fenn, Rourke and Sanderson [FRS-1, FRS-2] but taking into account the first Reidemeister move and degeneracy. We explain more in next subsections.

## 2.4. **Presimplicial sets and modules.**

Let $X_n$, $n \geq 0$ be a sequence of sets and $d_i = d_{i,n} : X_n \to X_{n-1}$, $0 \leq i \leq n$ maps (called face maps) such that:

(1) $d_i d_j = d_{j-1} d_i$ for any $i < j$.

Then the system $(X_n, d_i)$ satisfying the above equality is called a presimplicial set[14]. Similarly if $C_n$, $n \geq 0$ is a sequence of $k$-modules, for fixed commutative ring $k$ (e.g. $C_n = kX_n$) and $d_i = d_{i,n} : C_n \to C_{n-1}$, $0 \leq i \leq n$ are homomorphisms satisfying

(1) $d_i d_j = d_{j-1} d_i$ for any $i < j$,

then $(C_n, d_i)$ satisfying the above equality is called a presimplicial module. The important basic observation is that if $(C_n, d_i)$ is a presimplicial module then $(C_n, \partial_n)$, for $\partial_n = \sum_{i=0}^{n} (-1)^i d_i$, is a chain complex.

## 2.5. **One-term distributive homology.**

---

[14]The concept was introduced in 1950 by Eilenberg and Zilber under the name *semi-simplicial complex* [E-Z].

**Definition 2.2.** [Prz-6] *We define a (one-term) distributive chain complex* $\mathcal{C}^{(*)}$ *as follows:* $C_n = \mathbb{Z}X^{n+1}$ *and the boundary operation* $\partial_n^{(*)} : C_n \to C_{n-1}$ *is given by:*

$$\partial_n^{(*)}(x_0, ..., x_n) = (x_1, ..., x_n)+$$

$$\sum_{i=1}^{n}(-1)^i(x_0 * x_i, ..., x_{i-1} * x_i, x_{i+1}, ..., x_n).$$

The homology of this chain complex is called a one-term distributive homology of $(X; *)$ (denoted by $H_n^{(*)}(X)$).

We directly check that $\partial^{(*)}\partial^{(*)} = 0$, however it is useful to note that $(X^{n+1}, d_i)$ is a presimplicial (semi-simplicial) set with $d_i(x_0, ..., x_n) = (x_0 * x_i, ..., x_{i-1} * x_i, x_{i+1}, ..., x_n)$.

2.6. **Two-term rack (spindle) homology.** The trivial quandle $(X, *_0)$, is defined by $a *_0 b = a$. The 2-term rack homology of a spindle (right self-distributive system (RDS)) is defined by the presimplicial module $(C_n, d_i^R)$ with $d_i^R = d_i^{(*)} - d_i^{(*_0)}$. One can recognize here a precubic set $(X^n, d_i^\epsilon)$ (compare [Prz-7]). This precubic set and its geometric realization were important in the initial approach in [FRS-1, FRS-2].

2.7. **More general colorings.** We can consider more general colorings when we allow two parts of an overcrossing to have different colors. Namely, we start from the set of colors $X$ and consider a function $R: X \times X \to X \times X$ such that for any coloring of semi-arcs by colors from $X$ at any crossing the convention given in Figure 2.7 holds (here $R(a, b) = (R_1(a, b), R_2(a, b))$.

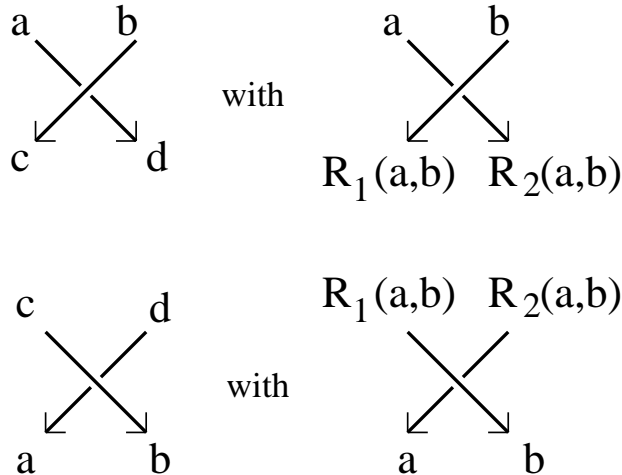

Figure 2.7: general semi-arc colorings

If $R$ is invertible and the number of colorings (for finite $X$) is preserved by a braid like oriented third Reidemeister move, we call $R$ a set theoretic Yang-Baxter operator which can be used to construct link invariants (e.g. 2-cocycle invariants) [CES, P-W-1].

Before we move to general Yang-Baxter operators and their invariants, I would suggest the reader the following simple but important exercise:

**Exercise 2.3.** *Let $D$ be an oriented link diagram, $X$ a finite set of colors, and $R : X \times X$ a set theoretic Yang-Baxter operator. Let $\phi$ be a coloring of semi-arcs of $D$ satisfying rules of Figure 2.7. Following Subsection 2.3, we define $\Psi(D, \phi) = \sum_v sgn(v)(a, b)$. Let us consider $\partial_2(a, b) = a + b - R_1(a, b) - R_2(a, b)$. Show that $\Psi(D, \phi)$ is a cycle, that is $\partial_2(\Psi(D, \phi)) = 0$ for every diagram $D$ and coloring $\phi$.*

## 3. Yang-Baxter homology

In quandle coloring and set-theoretic Yang-Baxter coloring of an oriented link diagram we are assuming that at every crossing a coloring of the input semi-arcs defines uniquely coloring of the output semi-arcs. We can however, in a natural way, relax this condition by allowing any coloring and then for a crossing to associate a weight from a fixed commutative ring (for set-theoretic Yang-Baxter operator this weight is 0 if coloring is not allowed and 1 if it is allowed). The details are as follows.

Fix a finite set $X$ and color semi-arcs of an oriented diagram $D$ by elements of $X$ allowing different weights from a fixed ring $k$ for every crossing. Following statistical mechanics terminology we call these weights Boltzmann weights. We allow also differentiating between a negative and a positive crossing; see Figure 2.8.

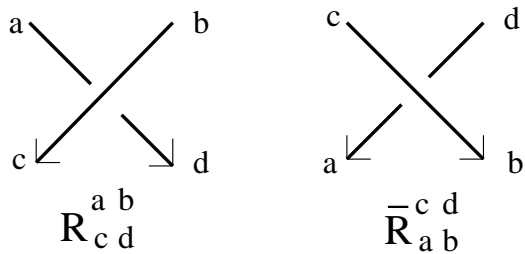

Figure 2.8; Boltzmann weights $R_{c,d}^{a,b}$ and $\bar{R}_{a,b}^{c,d}$ for positive and negative crossings

We can now generalize the number of colorings to state sum (basic notion of statistical physics) by multiplying Boltzmann weight over all crossings and adding over all colorings [Jon, Tur]:

$$col_{(X;BW)}(X) = \sum_{\phi \in col_X(D)} \prod_{p \in \{crossings\}} \hat{R}_{c,d}^{a,b}(p)$$

where $\hat{R}_{c,d}^{a,b}$ is $R_{c,d}^{a,b}$ or $\bar{R}_{c,d}^{a,b}$ depending on whether $p$ is a positive or negative crossing. Our state sum is an invariant of a diagram but to get a link invariant we should test it on Reidemeister moves. To get analogue of a shelf invariant we start from the third Reidemeister move with all positive crossings. Recall that in the distributive case, passing through a positive crossing was coded by a map $R : X \times X \to X \times X$ with $R(a, b) = (b, a * b)$. Thus in the general case passing through a positive crossing is coded by a linear map $R : kX \otimes kX \to kX \otimes kX$ and in basis $X$ the map $R$ is given by the $|X|^2 \times |X|^2$ matrix with entries $(R_{c,d}^{a,b})$, that is $R(a, b) = \sum_{(c,d)} R_{c,d}^{a,b} \cdot (c, d)$. The third Reidemeister move leads to the equality of the following maps $V \otimes V \otimes V \to V \otimes V \otimes V$ where $V = kX$:

$$(R \otimes Id)(Id \otimes R)(R \otimes Id) = (Id \otimes R)(R \otimes Id)Id \otimes R),$$

as illustrated in Figure 2.10. This is called the Yang-Baxter equation[15] and $R$ is called a pre-Yang-Baxter operator. If $R$ is additionally invertible it is called a Yang-Baxter operator. If entries of $R^{-1}$ are denoted by $\bar{R}_{c,d}^{a,b}$ then the state sum is invariant under "parallel" (directly oriented) second Reidemeister move, see Figure 2.9.[16]

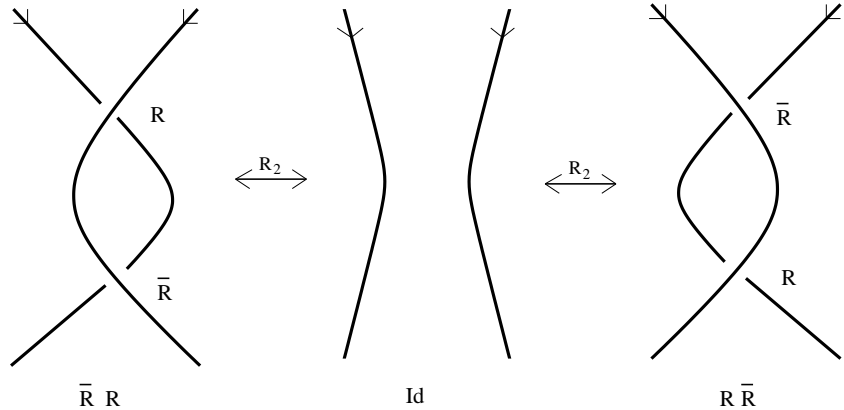

Figure 2.9; Invertibility of $R$ and the parallel second Reidemeister move

[16]We should stress that to find link invariants it suffices to use directly oriented second and third Reidemeister moves in addition to both first Reidemeister moves, as we can restrict ourselves to braids and use the Markov theorem. This point of view was used in [Tur].

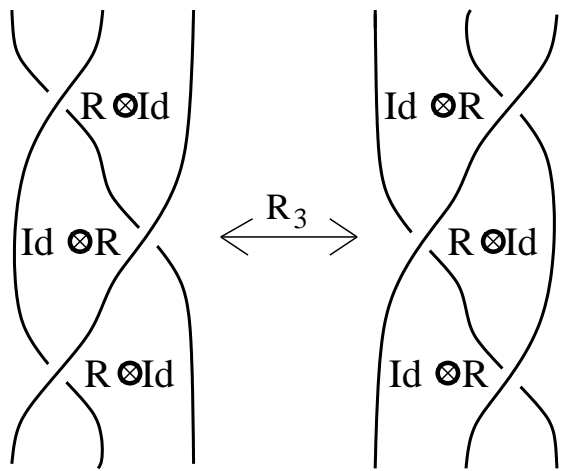

Figure 2.10; Yang-Baxter equation
from the positive third Reidemeister move

Examples leading to the Jones polynomial [Jon, Tur] start from a
2-dimensional set $X$ and the free $k$-module over $X$, $V = kX^2$ and
$R : V \otimes V \to V \otimes V$ is given by:

$$\begin{pmatrix} -q & 0 & 0 & 0 \\ 0 & q^{-1} - q & 1 & 0 \\ 0 & 1 & 0 & 0 \\ 0 & 0 & 0 & -q \end{pmatrix}$$

or using column unital (i.e. entries of each column adds to 1) matrix
[P-W-2, Wa]

$$\begin{pmatrix} 1 & 0 & 0 & 0 \\ 0 & 1 - y^2 & 1 & 0 \\ 0 & y^2 & 0 & 0 \\ 0 & 0 & 0 & 1 \end{pmatrix}$$

**Graphical visualization of Yang-Baxter face maps.**
The presimplicial set corresponding to a (two term) Yang-Baxter ho-
mology has the following visualization. In the case of a set-theoretic
Yang-Baxter equation we recover the homology of J. S. Carter, M. El-
hamdadi, M. Saito [CES]; compare [Leb-1, Leb-2, Prz-7, P-W-1].

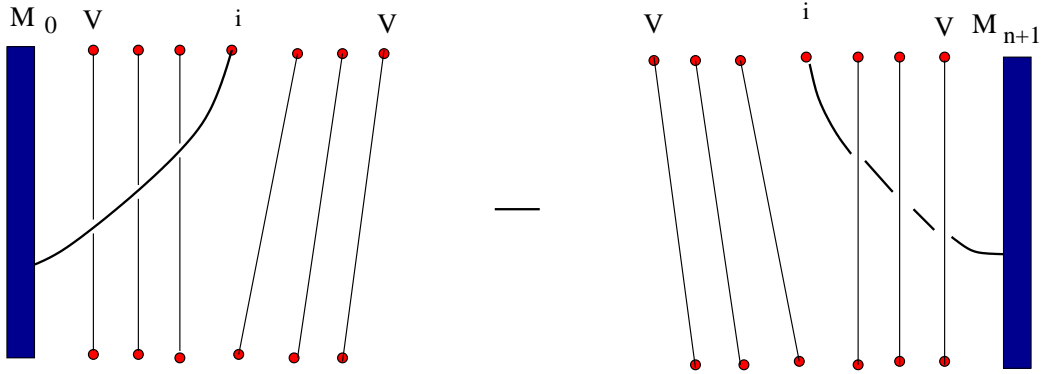

Diagramatic interpretation of a face map d $_i^{\text{Y-B}}$

Figure 3.1; Graphical interpretation of the face map $d_i$

In particular for a Yang-Baxter operator $R$ given by

$$R(a,b) = \sum_{(c,d) \in X^2} R_{c,d}^{a,b} \cdot (c,d)$$

we have

$$\partial_2(a,b) = (d_1^\ell(a,b) - d_1^r(a,b)) - (d_2^\ell(a,b) - d_2^r(a,b)) =$$

$$((b) - \sum_{c,d} R_{c,d}^{a,b} \cdot (c)) - (\sum_{c,d} R_{c,d}^{a,b} \cdot ((d) - (a)) =$$

$$(a) + (b) - \sum_{c,d} R_{c,d}^{a,b} \cdot ((c) + (d)).$$

**Exercise 3.1.** *Consider a Yang-Baxter operator $R$ and a coloring $\phi$ such that associated to a diagram elements $R_{c,d}^{a,b}$ are all different than zero. Find when the 2-chain $\Psi(D,\phi) = \sum_v sgn(v) \cdot (a,b)$ is a 2-cycle.*

## Decomposition of the third Reidemeister move into cubic face maps

The main idea is illustrated by the following picture, we can contemplate a precubic structure of the third Reidemeister move:

Yang–Baxter operator with $R = (R^{a,b}_{c,d})$ and fixed semi–arc coloring

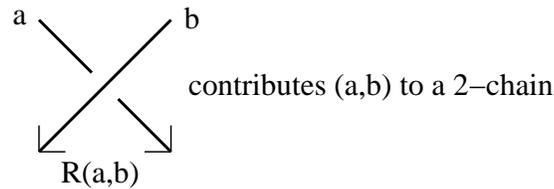

contributes (a,b) to a 2–chain

$$\partial_3(a,b,c) = \sum_{i=1}^{3} (-1)^i \, (d^\ell_i - d^r_i)(a,b,c)$$

We illustrate here the fact that the third Reidemeister move preserves homology classes, that is changes any chain by a boundary

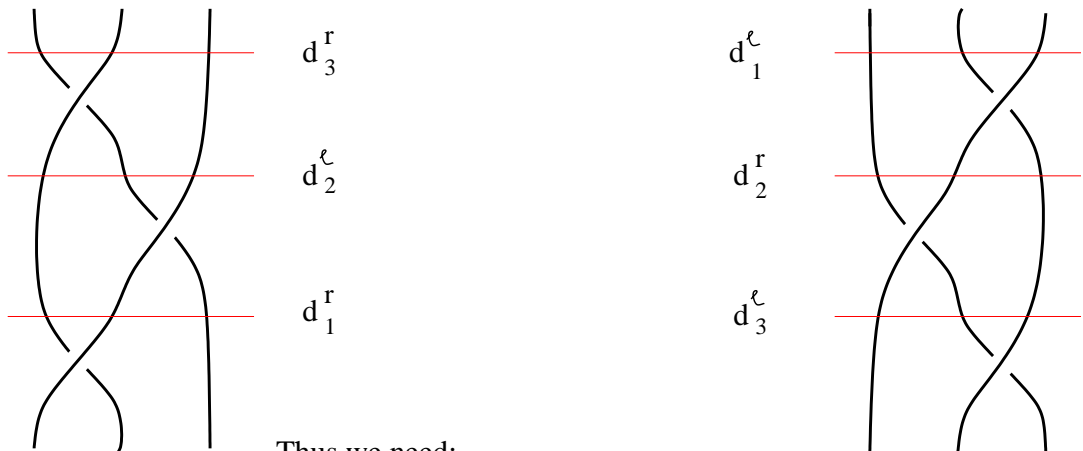

Thus we need:

$$d^r_3 + d^\ell_2 + d^r_1 = d^\ell_1 + d^r_2 + d^\ell_3$$

and this is given by $\partial_3(a,b,c)$

Figure 3.2; Reidemeister third move and face maps $d^\varepsilon_i$

The idea leads to (co)cycle invariants of links, at least for stochastic (or more generally column unital) Yang-Baxter matrices. An example was given in the Knot in Hellas talk by Xiao Wang (compare [P-W-2] and Wang's PhD thesis [Wa]).

## 4. Khovanov homology after Oleg Viro

One of the biggest discovery (or construction) in Topology after the first Knots in Hellas conference[17] was that of Khovanov homology [Kh].

We start from the description of the Khovanov homology for framed links, after [Vi-2, Vi-3].

**Definition 4.1.** *([Ka-1, Ka-2, Ka-3]) The unreduced Kauffman polynomial is defined by initial conditions*

$$[U_n] = (-A^2 - A^{-2})^n,$$

*where $U_n$ is the crossingless diagram of a trivial link of $n$ components, and the skein relation:*

$$[\,\times\,] = A[\,\asymp\,] + A^{-1}[\,)(\,].$$

A Kauffman state is a function from a set of crossings to the two element set $\{A, B\}$, that is $s : cr(D) \to \{A, B\}$, see Figure 4.1. We denote by $D_s$ the diagram (system of circles) obtained from $D$ by smoothing all crossings of $D$ according to $s$; $|D_s|$ denotes the number of circles in $D_s$.

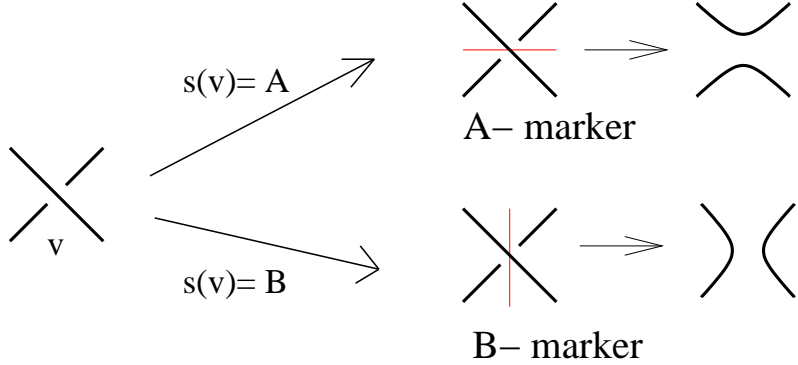

Figure 4.1, Interpretation of Kauffman states

**Proposition 4.2.** *(Kauffman) The unreduced Kauffman bracket polynomial can be written as the state sum (over all Kauffman states):*

$$[D] = \sum_{s \in KS} A^{|s^{-1}(A)| - |s^{-1}(B)|} (-A^2 - A^{-2})^{|D_s|},$$

*where $KS$ is the set of all Kauffman states of the diagram $D$.*

---

[17]The conference Knots in Hellas I took place in Delphi, Greece in August of 1998, while the e-print of Khovanov work was put on arXiv in August of 1999. However Mikhail Khovanov had already an idea of Khovanov homology in summer of 1997.

Notice that the Kauffman bracket associates to every trivial circle the polynomial $-(A^2 + A^{-2})$. In order to have state sum with monomial entries Viro considers two type of circles: positive with $A^2$ associated to it, and negative with $A^{-2}$ associated to it. These lead to Enhanced Kauffman States (EKS).

**Definition 4.3.**

(i) *An enhanced Kauffman state, $S$, is a Kauffman state $s$ together with a function $h : D_s \to \{+1, -1\}$.*

(ii) *The enhanced Kauffman state formula for the unreduced Kauffman bracket is the Kauffman state formula written using the set of enhanced Kauffman states EKS:*

$$[D] = \sum_{S \in EKS} (-1)^{|D_s|} A^{\sigma(s) + 2\tau(S)},$$

*where the signature of $s$ is $\sigma(s) = |s^{-1}(A)| - |s^{-1}(B)|$ and $\tau(S) = |h^{-1}(+1)| - |h^{-1}(-1)|$ that is the number of positive circles minus the number of negative circles in $D_s$ with enhanced Kauffman state function $h$ of $S$.*

With the above preparation we can define the Khovanov chain complex and Khovanov homology of a diagram.

**Definition 4.4.** *Consider bidegree on the Enhanced Kauffman States as follows:*

$$\mathcal{S}_{a,b} = \{S \in EKS \mid \sigma(s) = a, \quad \sigma(s) + 2\tau(S) = b\}$$

(a) *The chain groups are free abelian groups with basis $\mathcal{S}_{a,b}$, that is $C_{a,b}(D) = \mathbb{Z}\mathcal{S}_{a,b}$.*

(b) *Boundary maps are $\partial_{a,b} : C_{a,b}(D) \to C_{a-2,b}(D)$ given by the formula:*

$$\partial_{a,b}(S) = \sum_{S' \in \mathcal{S}_{a-2,b}} (-1)^{t(S,S')} [S; S'] S'$$

*where $[S, S']$ is 1 or 0 and it is 1 if and only if the following two conditions hold:*

(i) *$S$ and $S'$ differ at exactly one crossings, say $v$, at which $s(v) = A$ and $s'(v) = B$. In particular $\sigma(s') = \sigma(s) - 2$.*

(ii) *$\tau(S') = \tau(S) + 1$ and common circles to $D_s$ and $D_{s'}$ have the same sign. The possible signs of circles around the crossing $v$ is shown in Figure 4.2.*

*To define the sign $(-1)^{t(S,S')}$ we need to order crossings of $D$. Then $t(S, S')$ is equal to the number of crossings with label $A$ smaller than the crossing $v$ in the chosen ordering.*

(c) *The Khovanov homology is defined in the standard way as:*
$H_{a,b}(D) = ker(\partial_{a,b})/im(\partial_{a+2,b})$.

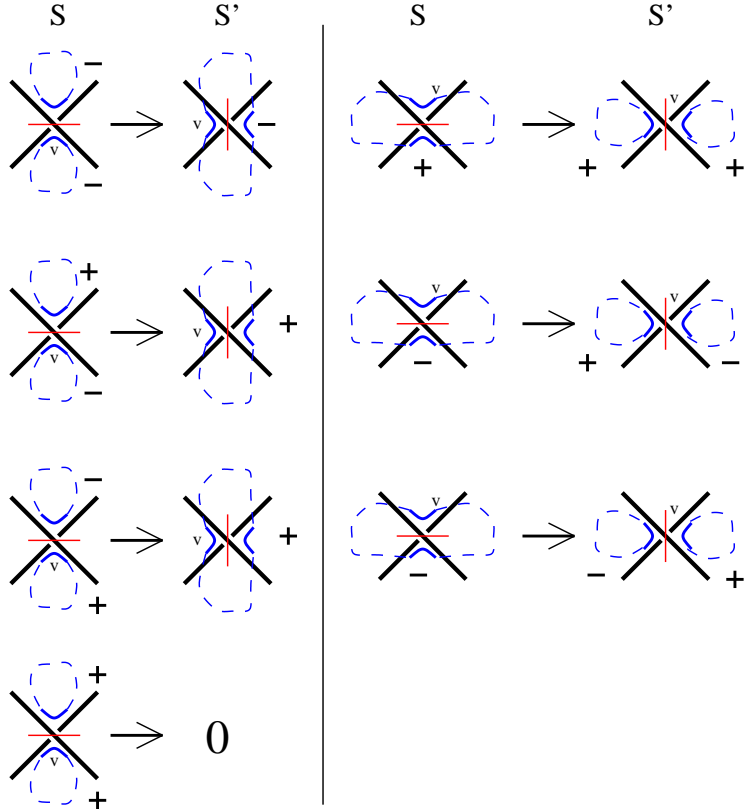

Figure 4.2, List of neighboring states with $[S, S'] = 1$

For every diagram $D$ we check easily that
$C_{a,b}(D) = 0$ for $a > cr(D)$ or $a < -cr(D)$, or $b > cr(D) + 2|D_{s_A}|$, or
$b < -cr(D) - 2|D_{s_B}|$.
These justify notation $a_{max} = cr(D)$, $b_{max} = cr(D) + 2|D_{s_A}|$, $a_{min} = -cr(D)$, and $b_{min} = -cr(D) - 2|D_{s_B}|$.
We always have $C_{a_{max},b_{max}} = \mathbb{Z} = C_{a_{min},b_{min}}$ but it often happens that $H_{*,b_{max}} = 0$ or $H_{*,b_{min}} = 0$. To say more we recall, after Lickorish and Thistlethwaite [L-T], the concept of adequate diagrams.

**Definition 4.5.** [L-T, AP] *We say, that a diagram $D$ is $s$-adequate for a Kauffman state $s$ if circles of $D_s$ have no self-touchings. Equivalently, $D$ is $s$-adequate if any diagram $D'_s$ obtained from $D$ by smoothing according to $s$ all but one crossing has smaller number of components than $D_s$. In particular, $D$ is said to be A-adequate if the state $s_A$ having all marker $A$ is adequate. Similarly, $D$ is said to be B-adequate if $s_B$ is adequate.*

We have classical observation ([Kh, Pr-Sa]):

**Proposition 4.6.** *For an A-adequate diagram D we have:*

$$H_{*,cr(D)+2|D_{s_A}|} = H_{cr(D),cr(D)+2|D_{s_A}|} = \mathbb{Z}$$

*Similarly for a B-adequate diagram D we have:*

$$H_{*,-cr(D)-2|D_{s_B}|} = H_{-cr(D),-cr(D)-2|D_{s_B}|} = \mathbb{Z}$$

In the case of any diagram $D$ the groups $H_{*,cr(D)+2|D_{s_A}|}$ and $H_{*,-cr(D)-2|D_{s_B}|}$ were studied in the PhD thesis of Marithania Silvero [S-C]. In particular she conjectures that these groups have no torsion. More on these groups and their geometric realization (conjectured to be of homotopy type of wedge of spheres) can be read in [GMS, Pr-Si]. We propose to call diagrams with nonzero Khovanov homology at $H_{*,cr(D)+2|D_{s_A}|}$ Khovanov $A$-adequate. Similarly Khovanov $B$-adequate diagram has to have nonzero groups $H_{-cr(D),-cr(D)-2|D_{s_B}|}$. We proved, playing odd Khovanov homology versus even Khovanov homology that there are links without $A$-Khovanov adequate diagrams. The simplest such example, we were able to find is the torus knot of type $T(4,-5)$.

Being in Greece let us look at some properties of Khovanov homology for the link of Lerna. Below is the table of its Khovanov homology computed by Sujoy Mukherjee using KhoHo program [Shu]. Since in tables one uses original oriented version of Khovanov (co)homology recall that if $\vec{D}$ is any oriented diagram of $D$ and $w(\vec{D})$ its writhe or Tait number then $H^{i,j}(\vec{D}) = H_{a,b}(D)$ for $i = \frac{w-a}{2}$ and $j = \frac{3w-b}{2}$.

The Lerna link has two components, and we can orient it so it has either $w(\vec{D}) = 16$ and all crossings positive (as in Figure 4.4) or $w(\vec{D}) = -4$. We use the second case in table calculation (Figure 4.3). With this convention we get the following unreduced Jones polynomial of the Lerna link[18]:

$$q^8 - 5q^6 + 12q^4 - 20q^2 + 27 - 29q^{-2} + 26q^{-4} - 18q^{-6} + 4q^{-8} + 11q^{-10} - 21q^{-12} +$$
$$29q^{-14} - 27q^{-16} + 23q^{-18} - 16q^{-20} + 10q^{-22} - 4q^{-24} + q^{-26}) =$$
$$(q + q^{-1})\Big( q^7 - 6q^5 + 18q^3 - 38q + 65q^{-1} - 94q^{-3} + 120q^{-5} - 138q^{-7} +$$
$$142q^{-9} - 131q^{-11} + 110q^{-13} - 81q^{-15} + 54q^{-17} - 31q^{-19} + 15q^{-21} - 5q^{-23} + q^{-25} \Big).$$

The second factor of the product is the reduced Jones polynomial of the Lerna link. Thus we notice that the coefficients alternate in signs, as the Lerna link is alternating. Furthermore the absolute values of the coefficients form a strictly unimodal sequence

$$1 < 6 < 18 < 38 < 65 < 94 < 120 < 138 < 142$$

$$142 > 131 > 110 > 81 > 54 > 31 > 15 > 5 > 1$$

which is also strictly logarithmically concave (i.e. $c_i^2 > c_{i-1}c_{i+1}$).

---

[18]To get the classical Jones notation we put $q = -t^{1/2}$.

| j \ i | -10 | -9 | -8 | -7 | -6 | -5 | -4 | -3 | -2 | -1 | 0 | 1 | 2 | 3 | 4 | 5 | 6 |
|---|---|---|---|---|---|---|---|---|---|---|---|---|---|---|---|---|---|
| 8 | | | | | | | | | | | | | | | | | 1 |
| 6 | | | | | | | | | | | | | | | | 5 | $1_2$ |
| 4 | | | | | | | | | | | | | | | 13 | $1,5_2$ | |
| 2 | | | | | | | | | | | | | | 25 | $5,13_2$ | | |
| 0 | | | | | | | | | | | | | 40 | $13,25_2$ | | | |
| -2 | | | | | | | | | | | | 54 | $25,40_2$ | | | | |
| -4 | | | | | | | | | | | 66 | $40,54_2$ | | | | | |
| -6 | | | | | | | | | | 73 | $55,65_2$ | | | | | | |
| -8 | | | | | | | | | 69 | $65,73_2$ | | | | | | | |
| -10 | | | | | | | | 62 | $73,69_2$ | | | | | | | | |
| -12 | | | | | | | 48 | $69,62_2$ | | | | | | | | | |
| -14 | | | | | | 33 | $62,48_2$ | | | | | | | | | | |
| -16 | | | | | 21 | $48,33_2$ | | | | | | | | | | | |
| -18 | | | | 10 | $33,21_2$ | | | | | | | | | | | | |
| -20 | | | 5 | $21,10_2$ | | | | | | | | | | | | | |
| -22 | | | $10,5_2$ | | | | | | | | | | | | | | |
| -24 | 1 | 5 | | | | | | | | | | | | | | | |
| -26 | 1 | | | | | | | | | | | | | | | | |

Figure 4.3; Table of Khovanov homology for the Lerna link

As the Lerna link is a non-split alternating link the nontrivial entries of Khovanov Homology are on two diagonals of slope 2 and torsion is on the lower one (Eun Soo Lee [Lee-1, Lee-2, Lee-3]. Furthermore we observed that there is only $\mathbb{Z}_2$ torsion. It agrees with a general, but yet not published, result of Alexander Shumakovitch that alternating links can have only $\mathbb{Z}_2$-torsion [Shu-1]. The adequacy of Lerna link is reflected in extreme coefficients $H^{-10,-26} = \mathbb{Z} = H^{6,8}$ if $w(\vec{D}) = -4$. Furthermore by results of [AP, PPS, Pr-Sa] and the fact that $D(Lerna)$ is strongly $A$-adequate we know that $H^{-10,-24} = \mathbb{Z}$ and $H^{-9,-24} = \mathbb{Z}^5$ (the Tait diagrams of the Lerna link are shown in Figure 4.4; notice that the A-smoothing diagram has no odd cycles).

$H^{-8,-22} = \mathbb{Z}^{10} \oplus \mathbb{Z}_2^5$ and $tor H^{5,4} = \mathbb{Z}_2^5$. To compute torsion here we use the following result from [Pr-Sa]. To formulate Theorem 4.8 we need to recall the notion of a state graph (Definition 2.1 of [AP]):

**Definition 4.7.** [AP, PPS] *Given a diagram $D$ and a Kauffman state $s$, we define an associated state graph $G_s$ with vertices in bijection with circles of $D_s$ and edges in bijection with crossings of $D$. An edge connects given vertices if the corresponding crossing connects circles of $D_s$ corresponding to the vertices.*

**Theorem 4.8.** [Pr-Sa] *For a given loopless graph $G$ let $G'$ denote the simple graph obtained from $G$ by replacing a multiple edge by a single one (see Figure 4.4).*

(A) *Let $D$ be an $A$-adequate diagram of $n$ crossings and $G_{s_A}$ associated graph Assume that $G_{s_A}$ is connected then:*

$$tor H_{n-4,n+2|D_{s_A}|-8} = \begin{cases} \mathbb{Z}_2^{p_1(G'_{s_A}(D))-1} & if \quad G'_{s_A}(D) \ has \ an \ odd \ cycle \\ \mathbb{Z}_2^{p_1(G'_{s_A}(D))} & if \quad G'_{s_A}(D) \ is \ a \ bipartite \ graph \end{cases}$$

(B) *Let $D$ be an $B$-adequate diagram of $n$ crossings and $G_{s_B}$ associated graph, then:*

$$tor H_{-n+2,-n-2|D_{s_B}|+8} = \begin{cases} \mathbb{Z}_2^{p_1(G'_{s_B}(D))-1} & if \quad G'_{s_B}(D) \ has \ an \ odd \ cycle \\ \mathbb{Z}_2^{p_1(G'_{s_B}(D))} & if \quad G'_{s_B}(D) \ is \ a \ bipartite \ graph \end{cases}$$

*Proof.* Part (A) is Proposition 4.8(i) of [Pr-Sa] and part (B) follows from (A) by Khovanov duality and universal coefficient theory. That is if $\bar{D}$ denotes the mirror image of $D$ then $H_{-a,-b}(D) = H^{a,b}(\bar{D}) = free(H_{a,b}(\bar{D})) \oplus tor(H_{a-2,b}(\bar{D}))$. $\square$

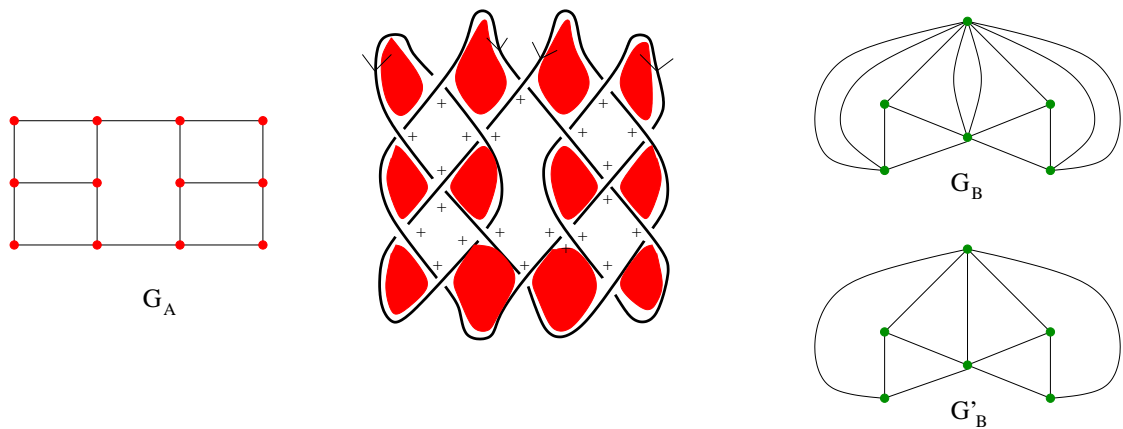

Figure 4.4, checkerboard coloring of Lerna diagram and Tait's graphs $G_{s_A}$ and $G_{s_B}$

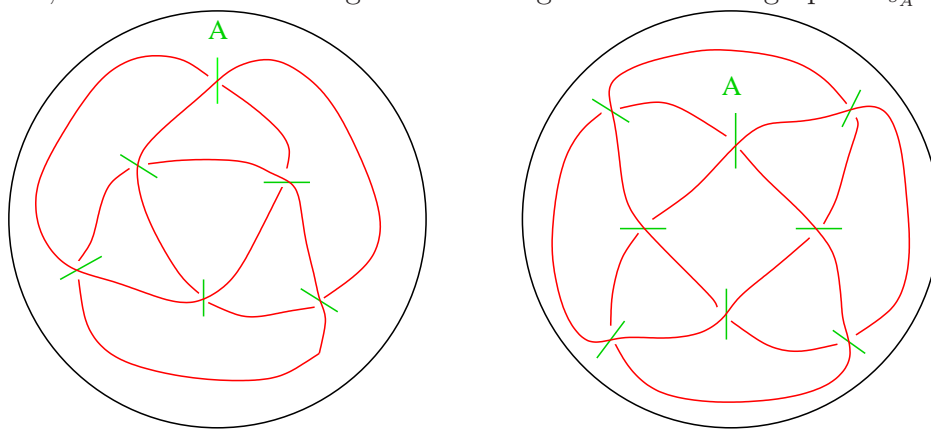

Figure 4.5; Links whose $A$ Kauffman states are Lerna pseudoknots

Notice that these links, closed 3-braids $(\sigma_1\sigma_2)^{-3}$ and $(\sigma_1\sigma_2)^{-4}$, respectively, are very far from being $A$-adequate but they are Khovanov $A$-adequate (see [S-C, GMS, Pr-Si]).

## 5. SUMMARY

We discussed historical ramifications of the beginning of Knot Theory. For the paper based on talks in Greece it is a natural turn. The reader can ask however what is a relation between distributive and Yang-Baxter homology on one hand and Khovanov homology on the other. The answer is simple: a connection is not yet found but my feeling is that we are only a step away. Maybe by the next Knots in Hellas III conference a link will be established and use of Khovanov homology in statistical physics will be demonstrated.

## 6. Acknowledgements

I would like to thank Sofia Lambropoulou for organizing for the second time the great Knots in Hellas conference.
I was partially supported by the Simons Collaboration Grant-316446 and CCAS Dean's Research Chair award.

Department of Mathematics,
The George Washington University,
Washington, DC 20052
e-mail: `przytyck@gwu.edu`,
and University of Gdańsk, Poland